# PREVALENCE: AN ADDENDUM

BRIAN R. HUNT, TIM SAUER, AND JAMES A. YORKE

Since the publication of our paper "Prevalence: a translation-invariant 'almost every' on infinite-dimensional spaces" in this journal [7], we have become aware of some work that is closely related to ours. We wish to call the reader's attention to this material. We thank J. Borwein, N. Kalton, and R. Dudley for informing us of this related work.

We defined the notions of "prevalent" and "shy" to be used in infinite-dimensional spaces of functions as analogues of the notions of "almost every" and "measure zero" with respect to Lebesgue measure on Euclidean spaces. Our definitions were given for complete metric linear spaces, and they extend trivially to abelian groups that are not vector spaces (but still have the topology of a complete metric). These definitions have been extended further to nonabelian groups by Mycielski [11].

For "abelian Polish groups" (topological abelian groups with a complete separable metric), Christensen [3] defined the notion of a "Haar zero set". Christensen's definition is equivalent to our definition of a shy set in the separable case; our definition has an extra provision which is only relevant for nonseparable spaces.

We believe there are many possible applications of these ideas. Christensen's main application [4, 5] is to prove an analogue of Rademacher's theorem (that a Lipschitz function from one Euclidean space to another is differentiable almost everywhere) for Lipschitz functions on Banach spaces. (Although this result is not true for the Fréchet derivative, it is for a slightly weaker notion of differentiability.) For other results concerning almost everywhere differentiability of Lipschitz functions on Banach spaces, see [1, 2, 9, 10, 12, 13]; some of these papers offer different notions of "almost everywhere". Christensen's definition has also been used [6, 8] in studying the continuity and differentiability of convex functions on Banach spaces. In [7] we presented ten results involving prevalence, with emphasis on dynamical systems and related areas such as transversality. The focus of our applications was different than that of the previous authors, who were primarily concerned with arbitrary Banach spaces. We were interested in proving that almost every function (or dynamical system) in a certain space (such as $C^1(\mathbb{R}^n)$) has a certain property. We made explicit the role that Lebesgue measure (on finite-dimensional subspaces) can play in proving such results, and we hope to have made accessible many more results of this type.

Finally, we would like to mention that Tsujii [14] has formulated a definition of "measure zero" for spaces of functions from one manifold to another. Tsujii gives several applications to transversality and dynamical systems [14–16].

(B. R. Hunt and J. A. Yorke) Institute for Physical Science and Technology, University of Maryland, College Park, Maryland 20742

*E-mail address*, B. Hunt: `hunt@ipst.umd.edu`

*E-mail address*, J. Yorke: `yorke@ipst.umd.edu`

(T. Sauer) Department of Mathematics, George Mason University, Fairfax, Virginia 22030

*E-mail address*: `tsauer@gmu.edu`

---

[1] The Ulam Quarterly is an electronic journal and is available in $\mathcal{AMS}$-TeX and PostScript format by anonymous ftp from `math.ufl.edu` in the directory `pub/ulam`. For information on ordering printed copies, contact: Professor Piotr Blass, Editor-in-Chief, Ulam Quarterly, P. O. Box 24708, W. Palm Beach, FL 33416-4708.